\documentclass[conference]{IEEEtran}
\IEEEoverridecommandlockouts
\usepackage{marvosym}

\usepackage[utf8]{inputenc} 
\usepackage[T1]{fontenc}    
\usepackage{hyperref}       
\usepackage{url}            
\usepackage{booktabs}       
\usepackage{amsfonts}       
\usepackage{nicefrac}       
\usepackage{microtype}      
\usepackage{lipsum}
\usepackage{fancyhdr}       
\usepackage{graphicx}       
\graphicspath{{media/}}     
\usepackage{dsfont}
\usepackage{amsmath}
\usepackage{eurosym}

\usepackage[normalem]{ulem}
\useunder{\uline}{\ul}{}

\usepackage{algorithm}
\usepackage{algpseudocode}

\usepackage{macros} 
\usepackage{tikz}

\DeclareMathOperator*{\argmin}{\text{arg}\,min}

\usepackage{xcolor}

\newcommand{\ie}{\textit{i.e. }}
 
\title{Controlling Microgrids Without External Data: \\ A Benchmark of Stochastic Programming Methods}

\author{\IEEEauthorblockN{Alban Puech$^1$, Tristan Rigaut$^1$, Adrien Le Franc$^2$, William Templier$^1$, Jean-Christophe Alais$^1$, Maud Tournoud$^1$,} \IEEEauthorblockN{ Victor Bossard$^1$, Alejandro Yousef$^1$, Elena Stolyarova$^1$}
\IEEEauthorblockA{\textit{$^1$ Schneider Digital - AI Hub, Schneider Electric} \\
\textit{$^2$ LAAS, CNRS}\\
France \\
\{alban.puech, tristan.rigaut, maud.tournoud\}@se.com}}

\begin{document}
\maketitle

\begin{abstract}
Microgrids are local energy systems that integrate energy production, demand, and storage units. They are generally connected to the regional grid to import electricity when local production and storage do not meet the demand. In this context, Energy Management Systems (EMS) are used to ensure the balance between supply and demand, while minimizing the electricity bill, or an environmental criterion. The main implementation challenges for an EMS come from the uncertainties in consumption, local renewable energy production, and in electricity price and carbon intensity. Model Predictive Control (MPC) is widely used to implement EMS but is particularly sensitive to the forecast quality, and often requires a subscription to expensive third-party forecast services. We introduce four Multistage Stochastic Control Algorithms relying only on historical data obtained from on-site measurements. We formulate them under the shared framework of  Multistage Stochastic Programming and benchmark them against two baselines in 61 different microgrid setups using the EMSx dataset \cite{le2021emsx}. Our most effective algorithm produces notable cost reductions compared to an MPC algorithm that utilizes the same uncertainty model to generate predictions, and it demonstrates similar performance levels to an ideal MPC that relies on perfect forecasts.  
\end{abstract}



\section{Introduction}

The US Department of Energy defines a microgrid as a ``group of interconnected loads and distributed energy resources with clearly defined boundaries that act as a single, controllable entity and can connect and disconnect from the grid to operate in both grid-connected and island mode'' \cite{TON201284}. Microgrids may play a prominent role in increasing climate resilience, with the opportunity to deploy more renewable energy and help to reduce energy costs.

Energy Management Systems (EMS) automatically control microgrids, optimizing an objective function (usually minimizing the cost), while respecting constraints under uncertainties in the energy production/demand/price. Model Predictive Control (MPC) is the most common method to implement an EMS \cite{mpc}, and simply computes controls (\ie actions to take) for a given time interval, using forecasts over a rolling horizon as proxies of future uncertainties. Since MPC is intuitive and scales well w.r.t the number of state variables, it has a wide history of industrial applications. Moreover, it relies on mature mathematical programming techniques and software. However, it is particularly sensitive to the quality of forecasts, so implementing a robust MPC for an industrial application can require paying third-party services to obtain them. These services are often expensive and increase the complexity of EMS development and maintenance.

Instead of being addressed by forecasts, uncertainties can be modelled as stochastic processes conditioned by past observations. 
This approach to EMS design has proved successful for various microgrid settings, e.g. in \cite{rigaut2018stochastic, le2021emsx, pacaud2022optimization}.
In all three references, stochastic programming methods \cite{bertsekas2012dynamic, shapiro2021lectures}
outperform MPC controllers in most test scenarios.
We find these results interesting from the perspective of removing third-party forecasts from the EMS components.

In this paper, we benchmark several EMS solutions
that are only allowed to make decisions based on past microgrid data observations.
We use a probabilistic modeling of future uncertainties
fitted on these data, which relies on the Darts toolbox \cite{darts}.
We compare four different algorithms, namely stochastic programming (\texttt{SP}), two-stage stochastic programming (\texttt{2S-SP}), two-stage stochastic programming with scenario clustering (\texttt{2S-SP-C}) and \texttt{MPC}. We further benchmark them against a heuristic (\texttt{HEU}) algorithm and an MPC algorithm relying on perfect forecasts (\texttt{P-MPC}). We assess them on a simple grid-connected microgrid consisting of a single energy storage, a local renewable production and a local load. Our main contributions are:
\begin{itemize}
    \item A comprehensive mathematical formulation of four different control methods that do not use any external data, under the shared framework of Multistage Stochastic Programming.
    \item A large-scale benchmark of the algorithms for the control of 61 anonymized industrial sites from the EMSx dataset made publicly available by Schneider Electric\footnote{\url{https://zenodo.org/record/5510400\#.YUizGls69hE}}~\cite{le2021emsx}.
    \item A quantitative analysis and discussion of the performance of the algorithms on the different sites of the dataset.
\end{itemize}

\label{sec:headings}

\section{Microgrid model and simulator}
\label{sec:simulator}


In this paper, we design an EMS to control a simple grid-connected microgrid. 
The microgrid includes a building consuming an uncertain load; solar panels producing an uncertain amount of energy; and a battery with an effective capacity $\overline{S} \in \mathbb{R}_+^*$, a maximum charge power $\overline{B} \in \mathbb{R}_+^*$ and a minimum discharge power $\underline{B} \in \mathbb{R}_-^*$. 

\subsection{Microgrid model}

We model the microgrid as a discrete-time dynamical system.
Decisions are made at every time-step of the sequence $\mathcal{T} = \{0,\ldots,T\}$
of length $T \in \mathbb{N}^*$ and time increment $\Delta t$.
We introduce the system dynamics during the time interval 
$\left[t, t+1 \right)$ of length $\Delta t$.

\paragraph{Battery operation}

 let $s_t$ denote the amount of energy stored in the battery at the beginning of the time interval; which is bounded by
 $0 \leq s_t \leq \overline{S}$. The EMS can charge/discharge an amount $p^b_t$ of energy in the battery, as expressed in the discrete-time dynamic equation:
\begin{equation} \label{eq:socdyn}
    s_{t+1} = s_t + \rho_c \max\{0,p^b_t\} + \rho_d^{-1} \min\{0,p^b_t\} 
    \eqfinp
\end{equation} 
$p^b_t$ may be zero, positive (charge), or negative (discharge) and is bounded by $\underline{B}\Delta t \leq p^b_t\leq \overline{B}\Delta t$.
 Each microgrid in the benchmark has specific
charge/discharge efficiency coefficients
$\rho_c \in (0,1]$ and $\rho_d \in (0,1]$, and an initial energy stock $s_0$.
 
\paragraph{Energy production and demand}
we denote by $p_t$ the energy produced by solar panels
and by $d_t$ the energy demand.

\paragraph{Energy exchanges with the grid} the EMS can import an energy amount 
\begin{equation}
    \label{eq:imports}%
    e_t = \max\{0,d_t - p_t + p^b_t\}
\end{equation}
from the power grid.
These imports occur when the \emph{net load} $d_t - p_t$ plus energy in/outflows $p^b_t$
from the battery are non-negative.

\paragraph{Grid import fee} electricity imports 
induce a cost that equates to the price of the energy volume $e_t$,
paid at a cost $c_t$ in \euro/kWh,
plus a fixed penalty of $\overline{C}$ in \euro\ if this volume exceeds the subscribed power $\overline{E}$ of the installation. It is expressed as:
\begin{equation}
    c_t \times e_t + \overline{C} \times \delta_{\geq \overline{E}}(e_t)
    \eqfinv
    \label{eq:bill}%
\end{equation}
where $\delta_{\geq \overline{E}}$ is the indicator function equal to $1$ when the import exceeds the maximum power.
We assume that when the microgrid produces more energy than it consumes, the exported energy is not remunerated by the grid operator.

\subsection{Optimal control framework}

\paragraph{Standard control notations}

we adopt standard control notations as used in \cite{bertsekas2012dynamic}.
State, control, and exogenous random variables are respectively noted as:
\begin{equation}
    x_t = s_t \eqsepv   
    u_t = p^b_t \eqsepv
    w_t = d_t - p_t \eqfinp
\end{equation}
The aforementioned dynamics~\eqref{eq:socdyn} and the constraints
on $(s_t, p^b_t)$ are summarized as follows:
\begin{equation}
    x_{t+1} = f(x_t, u_t, w_t) \eqsepv
    g_t(x_t,u_t,w_t) \leq 0 \eqfinv
\end{equation}
where $f_t: \RR^3 \to \RR_+$ and $g_t: \RR^3 \to \RR^{n_c}$ denote respectively the dynamics function and the constraints. 
The stage cost described in~\eqref{eq:bill}
is written as a function $L_t : \RR^3 \to \RR_+$
of $(x_t, u_t, w_t)$.

\paragraph{Controller definition}

At any time-step $t \in \mathcal{T}$,
the decision $u_t$ is made
based on the current state $x_t$, and on the \textit{history} of the microgrid
up to time $t$, defined as:
\begin{equation}
	h_t = \na{w_s}_{s \leq t}  \eqfinp
	\label{eq:history}%
\end{equation}
Then, a microgrid controller is a collection
$\phi = \na{\phi_t}_{t \in \Times}$
of mappings 
\begin{equation}
	\phi_t : (x_t, h_t) \mapsto u_t \eqsepv
	\forall t \in \Times \eqfinp
	\label{eq:controler}
\end{equation}
By design, this generic definition of a microgrid controller is agnostic to the method employed to control the system. 
This opens our benchmark to any other method from all fields of 
discrete-time control theory. The core distinction with the microgrid controller from the EMSx benchmark~\cite{le2021emsx} is that
the definition of history in equation~\eqref{eq:history}
does not include external forecasting data.

\paragraph{Evaluating the controllers}
We evaluate the performance of a controller by measuring the costs under its policy $\phi$ along the scenario 
$h = \na{h_t}_{t \in \Times}$, starting from an initial
given state $x_0$:

\begin{subequations}
\label{eq:cost}
    \begin{align}
    J(\phi, h) &= \sum_{t \in \Times} L_t(x_t, u_t, w_t) \eqfinv \\
    &
    \text{ where }
    \begin{cases}
    x_{t+1} = f(x_t, u_t, w_t) \eqsepv \\
    g_t(x_t,u_t,w_t) \leq 0 \eqfinv \\
    u_t = \phi_t(x_t, h_t) \eqfinv
    \end{cases}
    \forall t \in \Times \eqfinp
\end{align}    
\end{subequations}

\section{Presentation of the challengers}

In this section, we introduce our challengers: the different control methods that we have implemented to benchmark against a perfect forecast MPC and a heuristic controller. We frame them all as scenario tree-based Multistage Stochastic Programming (MSSP) policies that differ by the scenario tree that they use.





\subsection{Multistage Stochastic Programming policies}
At every time-step $t_0$ an MSSP policy observes the state $x_{t_0}$ and its history $h_{t_0}$ and searches a control $\phi_{t_0} (x_{t_0}, h_{t_0})$ in the argmin of the following stochastic optimization problem: 
\begin{subequations}\label{eq:generic-sp}
\begin{align}
    \argmin_{u_{t_0}}&  \min_{(\mathbf{U}_{t_0+i})_{1\leq i\leq R}} \mathbb{E} \sum_{t=t_0}^{t_0+R} L_{t}( \mathbf{X}_{t}, \mathbf{U}_{t}, \mathbf{W}_{t}) \eqfinv \label{eq:generic-sp-cost}\\
    \text{s.t. }& \mathbf{X}_{t+1} = f_{t}(\mathbf{X}_{t}, \mathbf{U}_{t}, \mathbf{W}_{t}) \eqfinv \label{eq:generic-sp-dyn} \\
    & g_{t}(\mathbf{X}_{t}, \mathbf{U}_{t}, \mathbf{W}_{t}) \leq 0\eqfinv\\
    & \mathbf{X}_{t_0}=x_{t_0}, \mathbf{H}_{t_0}=h_{t_0}\eqfinv \\
    & \sigma(\mathbf{U}_t) \subset \sigma(\mathbf{H}_{t_0},\mathbf{W}_{t_0+1},\ldots,\mathbf{W}_{t}) \eqfinp \label{eq:nonant}
\end{align}
\end{subequations}
where $R$ is a rolling horizon and the expectation integrates over the conditioned distribution $\mathbb{P}(\mathbf{W}_{t_0+1},\ldots,\mathbf{W}_{t_0+R} | h_{t_0})$. The policy only returns a control $u_{t_0}$ for time-step $t_0$, that is an element of the $\argmin$ in~\eqref{eq:generic-sp-cost}. The $\min$ over the subsequent decisions in~\eqref{eq:generic-sp-cost} measures the expected impact of our current control on the future costs incurred by these decisions. The measurability constraint~\eqref{eq:nonant} states these future decisions are functions of future uncertainties. If the distribution $\mathbb{P}(\mathbf{W}_{t_0+1},\ldots,\mathbf{W}_{t_0+R} | h_{t_0})$ has infinite support, the problem is hard to solve because the solution space is infinite-dimensional.

We introduce the concept of scenario tree and how it is used to approximate~\eqref{eq:generic-sp} by a tractable optimization problem.

\subsubsection{Scenario tree} \label{subsubsec:scenario-trees}

A scenario tree $G=(\mathcal{N}, \mathcal{P})$ is a tree graph representing the distribution of a discrete stochastic process $(\mathbf{\Xi}_0,\ldots,\mathbf{\Xi}_R)$ with finite support. $\mathcal{N}$ is a finite set of nodes, a node $n$ at a depth $t$ is a realization $(\xi_1,\ldots,\xi_{t})$ of the process up to time $t$. $\mathbf{\Xi}_0$ has only $1$ realization $\xi_0$ with probability $1$ so there is a single node, the root $n_0$, at depth $0$. Elements $(m,n) \in \mathcal{P} \subset \mathcal{N} \times \mathcal{N}$ are edges of the tree. We call $d_n$ the depth of a given node $n$. We define the probability $\pi(n)$ of a node $n$ as $\mathbb{P}(\mathbf{\Xi}_0=\xi_0,\ldots,\mathbf{\Xi}_{d_n}=\xi_{d_n})$, or equivalently, the product of probability transitions from the root $n_0$ to $n$. 
	
	
	

\subsubsection{Construction of a scenario tree} Many methods exist to generate a scenario tree modeling the distribution of a given stochastic process $(\mathbf{W}_0,\ldots,\mathbf{W}_R)$ from a collection of samples. We refer to~\cite{kammammettu2023scenario} and~\cite{kirui2020new} for recent reviews of the methods.

In this paper, we generate scenario trees to approximate the distribution of the netload stochastic process $(\mathbf{W}_{t+1},\ldots,\mathbf{W}_{t+R})$ conditioned by the netload history $h_t = (w_0,\ldots,w_t)$. To do so, we train a machine learning model to sample from $\mathbb{P}(\mathbf{W}_{t+1},\ldots,\mathbf{W}_{t+R} | h_t)$. Our model uses exclusively past realizations to model future uncertainties: this is an \emph{autoregressive} model. 

\subsubsection{Scenario tree approximation}

Approximating the stochastic process $\mathbb{P}(\mathbf{W}_{t_0+1},\ldots,\mathbf{W}_{t_0+R} | h_{t_0})$ by a scenario tree $G$ in \eqref{eq:generic-sp} turns the solution space into a finite dimensional one. Indeed constraint~\eqref{eq:nonant} implies that, for each node $n$ of $G$, there is one decision $u_n$, and due to constraint~\eqref{eq:generic-sp-dyn} one state $x_n$. We call $w_n = w_{d_n}$ the last netload realization of the node $n$. We rewrite~\eqref{eq:generic-sp} into its extensive formulation~\eqref{eq:sp}.
\begin{subequations} \label{eq:sp}
\begin{align}
    \argmin_{u_{n_0}}  &\min_{u_{n \in \mathcal{N} \backslash \{n_0\} }}  \sum_{n \in \mathcal{N}} \pi(n) \times L_{d_n+t_0}(x_{n}, u_{n}, w_{n}) \eqfinv \\
    \text{s.t. }& x_{n} = f_{d_m+t_0}(x_{m}, u_{m}, w_{m}),~\forall (m, n) \in \mathcal{P} \eqfinv \\
    &g_{d_n+t_0}(x_{n}, u_{n}, w_{n})\leq 0,~\forall n \in \mathcal{N} \eqfinv\\
    & x_{n_0} = x_{t_0} \eqfinp
\end{align}
\end{subequations}
In our case, problem~\eqref{eq:sp} is a 
so-called Mixed-Integer Linear Program (MILP) that we model with Pyomo \cite{PYOMO} and solve using HiGHS \cite{HIGHS}.
 The number of variables (resp. constraints) in~\eqref{eq:sp} grows linearly with the number of nodes (resp. edges) in the tree, while the latter grows exponentially with the number of time steps in the rolling horizon. This significantly deteriorates the performance of a MILP solver. 
 
 We present hereunder four methods to compute controls using scenario trees that are small enough to solve problem~\eqref{eq:sp} in a reasonable amount of time.

\subsection{Challengers}
\begin{figure}[!htbp]
\centering
\includegraphics[width=1.\columnwidth]{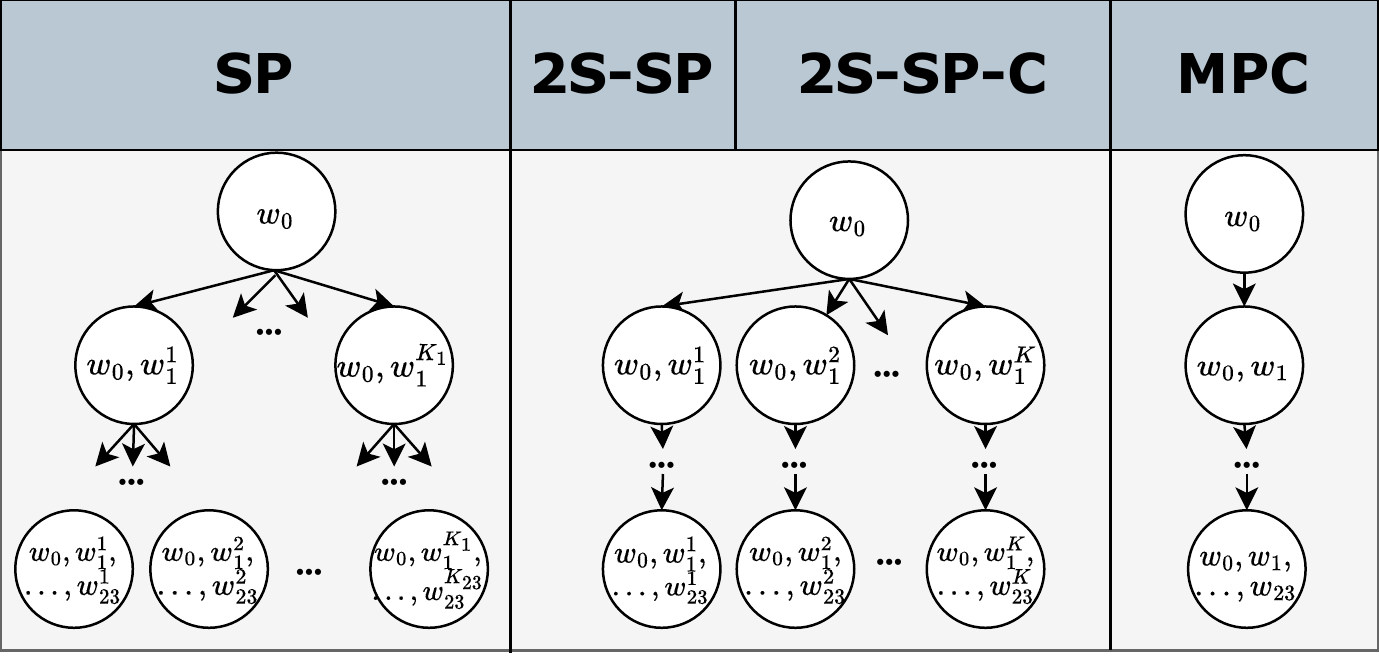}
\caption{Diagram of the construction of the trees for each algorithm. We use $R=23$ (control horizon of 24h)}
\label{fig:tree}
\end{figure}

Without loss of generality, we fix $R=23$ (control horizon of 24 hours).
We formulate each of our challengers as an MSSP method and use a LightGBM model implemented using Darts \cite{darts} to draw $K$ samples $(\tilde{w}_0,\tilde{w}_1^i,\tilde{w}_2^i,...,\tilde{w}_{23}^i),1\leq i \leq K$ from $\mathbb{P}(\mathbf{W}_{t+0},\ldots,\mathbf{W}_{t+R} | h_t)$, that are later used to generate a scenario tree. A challenger is identified by the kind of scenario tree it generates, as shown in Fig. \ref{fig:tree}:

\subsubsection{\texttt{SP}}

\texttt{SP} (short for Stochastic Programming) takes as input $(\tilde{w}_0,\tilde{w}_1^i,\tilde{w}_2^i,...,\tilde{w}_{23}^i),1\leq i \leq K$ and builds a scenario tree with $K_{23}$ leaves $(w_0,w_1^{1 \leq i_1\leq K_1},w_1^{1\leq i_{2}<K_{2}},...,w_{23}^{1\leq i_{23}<K_{23}}) $ using scenred  \cite{scenred}. The reduction parameter defines the number of nodes $K_1, K_2, ..., K_{23}$ at each depth of the tree. 

\subsubsection{\texttt{2S-SP} and \texttt{2S-SP-C}}

Two-stage stochastic programming, denoted by \texttt{2S-SP} directly turns each sample into a branch of the tree connected to the root so that $\forall 1\leq i \leq K:(w_0^i,w_1^i,...,w_{23}^i)=(\tilde{w}_0^i, \tilde{w}_1^i, ..., \tilde{w}_{23}^i)$. This type of tree is commonly referred to as a scenario fan. \texttt{2S-SP-C} is a slightly different version of \texttt{2S-SP} that uses K-means to separate the K samples into K' clusters, whose centres are used instead of the original samples to generate the tree.

\subsubsection{\texttt{MPC}}

Finally, Model predictive Control \texttt{MPC} is framed as an MSSP method, where the unique scenario (unique branch of the tree) is obtained by averaging all the samples: $(w_t)_{t\leq23} = \big((1/K) \sum_{1 \leq i \leq K} \tilde{w}_t^i\big)_{t\leq23} $

\subsection{Two baselines} 
The challengers are compared with two baseline controllers. 

\subsubsection{\texttt{P-MPC}} An ideal MPC with perfect forecasts that permits us to assess the value of perfect forecasts compared to our autoregressive models of uncertainties.

\subsubsection{\texttt{HEU}}
A myopic heuristic that stores electricity when the netload is negative and discharges the battery when the netload is positive. This baseline is called \texttt{HEU} and represents the traditional optimization-free industrial method.

\section{Results and discussions}

\begin{table}[!htbp]
\centering
\caption{Hyperparameter values and implementation specifics\label{table:hpp}}
\scalebox{0.9}{
\begin{tabular}{|ll|}
\hline
\multicolumn{2}{|c|}{{\ul \textbf{All challengers (MPC, SP, 2S, 2S-SP, 2S-SP-C)}}} \\ \hline
\multicolumn{1}{|l|}{\# past netload values used for the forecast (time-steps)}     & 48     \\ \hline
\multicolumn{1}{|l|}{Forecasting model}                                      & LightGBM               \\ \hline
\multicolumn{2}{|c|}{{\ul \textbf{SP}}}                                                               \\ \hline
\multicolumn{1}{|l|}{\# samples used for tree contruction} & 50                     \\ \hline
\multicolumn{1}{|l|}{Relative distance during (construction, reduction)}     & 0.2, 0.0                \\ \hline
\multicolumn{2}{|c|}{{\ul \textbf{2S-SP}}}                                                            \\ \hline
\multicolumn{1}{|l|}{\# samples used for tree contruction}                                   & 20                     \\ \hline
\multicolumn{2}{|c|}{{\ul \textbf{2S-SP-C}}}                                                          \\ \hline
\multicolumn{1}{|l|}{\# samples before clustering}                                   & 100                    \\ \hline
\multicolumn{1}{|l|}{\# clusters used for tree contruction}                          & $\leq$ 20 \\ \hline
\end{tabular}}
\end{table}

\subsection{Simulation Setup}
Our challengers (\texttt{MPC}, \texttt{SP}, \texttt{2S-SP} and \texttt{2S-SP-C}) are tested against the two baselines (\texttt{HEU} and \texttt{P-MPC}) through 61 simulations. Each simulation is associated with a time series of solar production and energy demand from a given site\footnote{Nine sites from \cite{le2021emsx} were not used for our simulations, as they contain extended periods of missing data}, with a 15-minute resolution, resampled into series with 1-hour-long time-steps. The energy demand time series corresponds to our exogenous uncertainty. Sites each have a microgrid setup, and we further set the subscribed power limit as  $\overline{E} = \max_{t}({d_t - p_t})-\overline{B}$. The penalty for exceeding the power subscription is fixed at 14.31\euro/h, as per the "Tarif Jaune" rates from France's main electric utility company EDF \cite{EDF}. The electricity prices (also displayed in Fig. \ref{fig:trajectories}) are set to 0.102\euro ~during off-peak hours (i.e. from 00:00 to 06:00, 09:00 to 11:00, 13:00 to 17:00, and 21:00 to 00:00) and to 0.153\euro ~during peak hours. We use 60 \% of the dataset for training, and the rest for testing. Table \ref{table:hpp} shows the hyperparameter values and other implementation specifics used for the different algorithms across the 61 simulations. Our control horizon is 24h ($R=23$).
All experiments were run on a computer with an Intel i7-12800H CPU and 32 GB RAM.

\subsection{Results on the 61 sites}

We introduce the average performances, both in terms of savings and computation times, of the challengers and baselines across the whole dataset of $61$ sites.

\begin{table*}[!htbp]
\caption{Simulation results}
\label{table:time}
\centering
\scalebox{0.85}{
\begin{tabular}{|c|c|c|c|c|c|c|}
\hline
\textbf{Alg.} & \begin{tabular}[x]{@{}c@{}}\textbf{Mean}\\ \textbf{processing}\\ \textbf{time (ms/it)}\end{tabular}  & \begin{tabular}[x]{@{}c@{}}\textbf{Avg}\\ \textbf{savings w.r.t} \\ \textbf{HEU (\%)}\end{tabular}  & \begin{tabular}[x]{@{}c@{}}\textbf{Extra}\\ \textbf{cost w.r.t} \\ \textbf{P-MPC (\%)}\end{tabular} & \begin{tabular}[x]{@{}c@{}}\textbf{Avg yearly} \\ \textbf{savings w.r.t} \\ \textbf{HEU (EUR)} \end{tabular} & \begin{tabular}[x]{@{}c@{}}\textbf{Times ranked} \\ \textbf{best challenger} \\ \textbf{algorithm (\%)}\end{tabular}  \\ \hline
\textbf{\texttt{MPC} }               & 36.16                                               & 2.77  & 3.86 & 13599   & 0                    \\ \hline
\textbf{\texttt{SP}}                 & 333.72                                               &\textbf{4.98}  & \textbf{1.51} & \textbf{21012} & \textbf{38}                      \\ \hline
\textbf{\texttt{2S-SP}}              & 249.35                                             & 4.94  & 1.56 & 20653 & 33                       \\ \hline
\textbf{\texttt{2S-SP-C}}            & 344.36                                             & 4.93  & 
1.57        & 20584    & 29          \\ \hline \hline
\textbf{\texttt{HEU}}                & 0.04                                               & .  & 6.87  & . & .                          \\ \hline
\textbf{\texttt{P-MPC}}              & 7.01                                               & 6.39  & . & 25945 & .                      \\ \hline

\end{tabular}}
\end{table*}

\begin{figure}[!htbp]
\centering
\includegraphics[width=0.9\columnwidth]{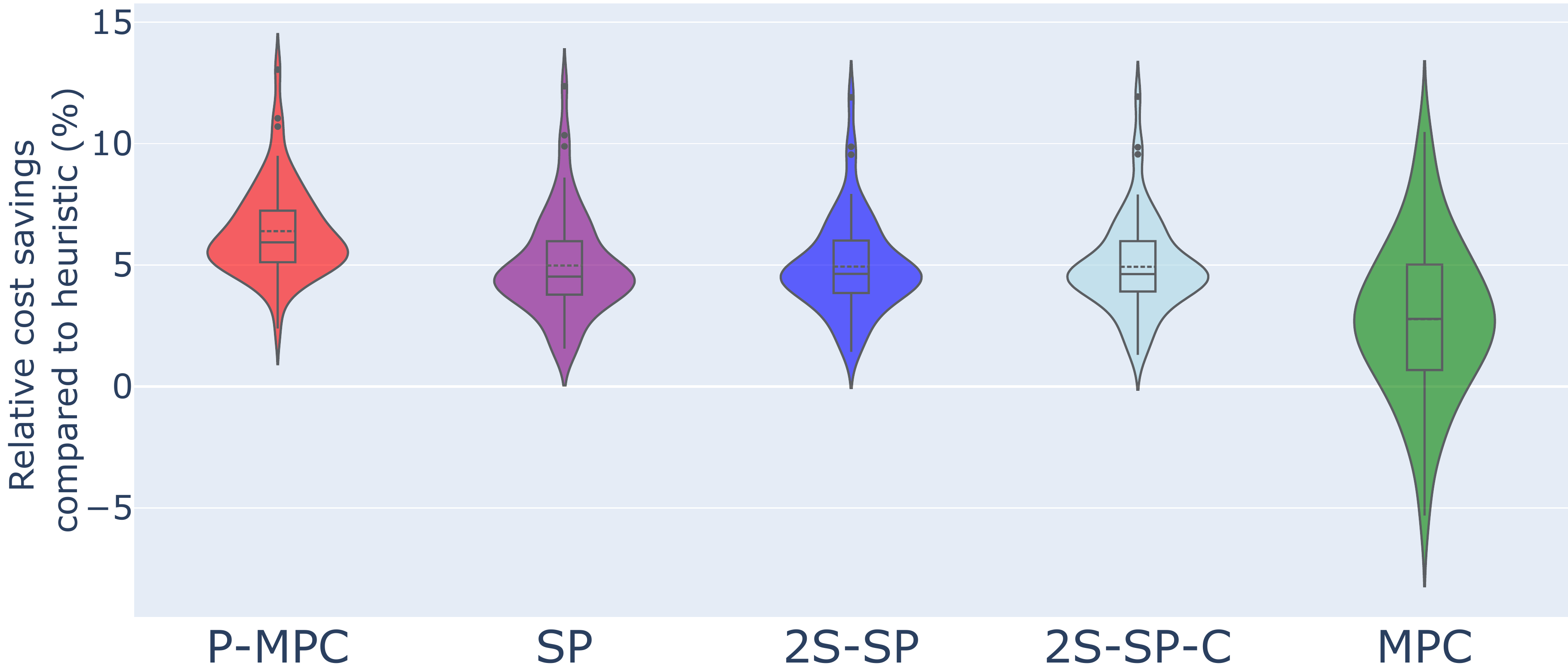}
\caption{Relative cost savings from using each of the challenger (\texttt{SP}, \texttt{2S-SP}, \texttt{2S-SP-C}) and \texttt{P-MPC} over the heuristic algorithm \texttt{HEU}}
\label{fig:rel_diff}
\end{figure}

\subsubsection{Average savings}

We show in Fig. \ref{fig:rel_diff} the distribution of the relative cost savings obtained using a given microgrid controller vs \texttt{HEU}. Using \texttt{SP}, \texttt{2S-SP}, and \texttt{2S-SP-C} over \texttt{HEU} always yields cost savings, and all three algorithms show similar performances. The savings obtained using  \texttt{MPC} are smaller on average, sometimes negative, and show a much larger variance. Note that, in contrast to the other challengers,  \texttt{MPC} is performing worse than \texttt{HEU} on almost 25 \% of the 61 different simulations.
We show in Table \ref{table:time} the average results obtained by the different challengers across the 61 simulations. \texttt{2S-SP} and \texttt{2S-SP-C} together outperform \texttt{SP} in the majority of the scenarios (71 \%), while \texttt{SP} yields higher cost savings on average (4.98 \% against 4.94 and 4.93 \% for \texttt{2S-SP} and \texttt{2S-SP-C} respectively). \texttt{MPC} sometimes performs worse than \texttt{HEU}, and the difference in performance between \texttt{P-MPC} and \texttt{MPC} highlights the importance of having accurate forecasts.  
Finally, the difference between \texttt{SP} and \texttt{P-MPC} was quite limited (1.51\%): using a probabilistic model only with local historical data alleviates the impact of an imperfect forecast.

\subsubsection{Computing times}

We additionally report in Table~\ref{table:time} the mean processing time per control. It corresponds to an average over a 7700-time-step-long test dataset. The training is done on 60 \% of the whole dataset and took 2.8s for each challenger.\texttt{HEU} computes a control in the smallest amount of time. \texttt{P-MPC} requires a bit more time as it creates and solves a MILP problem. \texttt{MPC} additionally uses the LightGBM model to compute forecasts. \texttt{2S-SP} also draws samples from the LightGBM model but creates a larger MILP problem than \texttt{MPC}, there are more nodes in the scenario tree. \texttt{2S-SP-C} adds a clustering step to \texttt{2S-SP} to solve a smaller MILP but still requires more time than \texttt{2S-SP}. Finally, \texttt{SP} stands between \texttt{2S-SP} and \texttt{2S-SP-C}, it adds the scenario tree generation to \texttt{2S-SP}, which takes less time than the clustering of \texttt{2S-SP-C}.

\subsection{Comparison of the control trajectories on a simulation sample}
\begin{figure}[!htbp]
\centering
\includegraphics[width=\columnwidth]{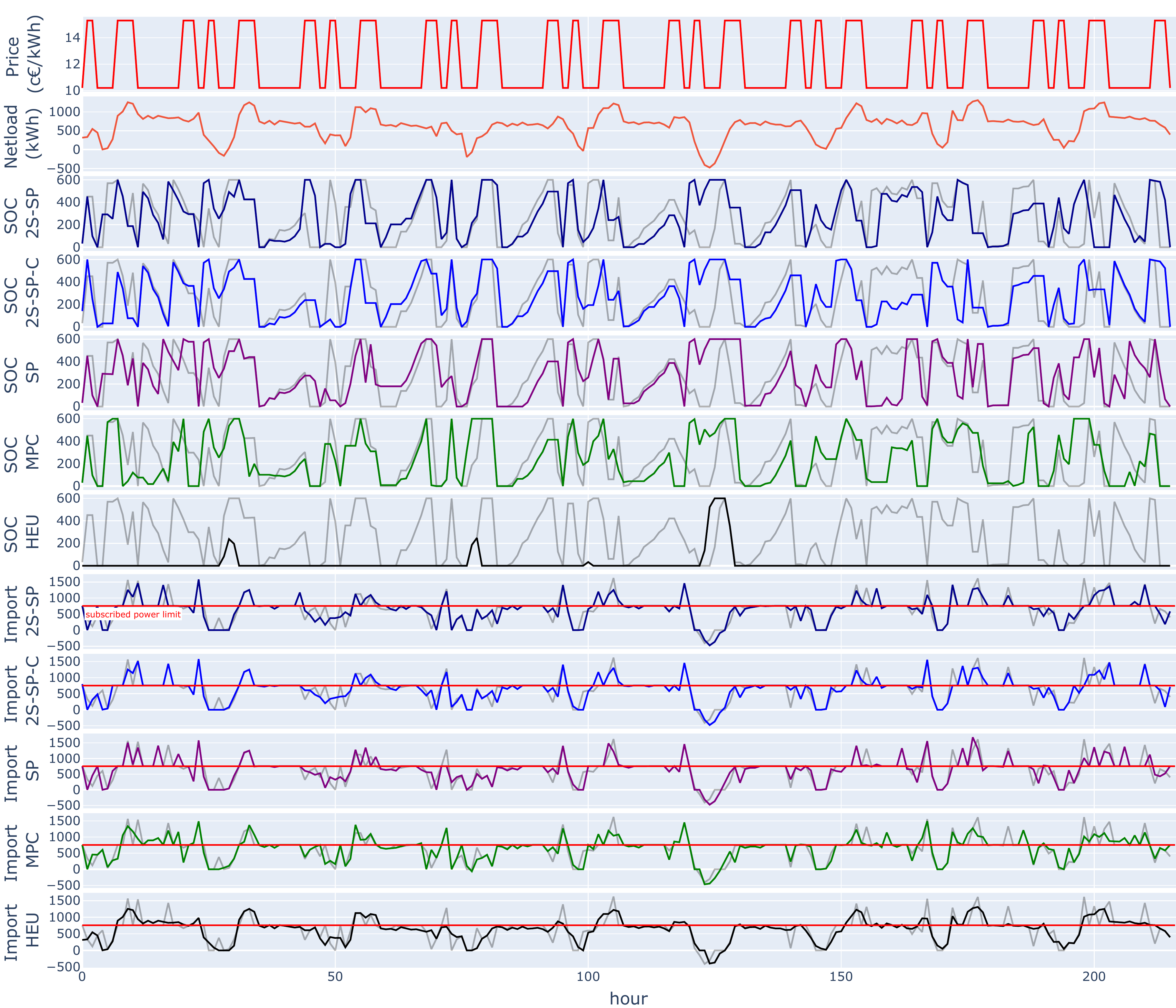}
\caption{Sample the control sequences/ trajectories of the challengers and the heuristic algorithm (\texttt{HEU}) running on site 19 of the EMSx dataset. "Import" refers to the amount of electricity imported from the grid, expressed in KWh.}
\label{fig:trajectories}
\end{figure}

We show in Fig. \ref{fig:trajectories} a sample of the battery state of charge and grid import trajectories under the policies of the different algorithms \footnote{Note that only a small portion of the simulation on site 19 of the EMSx dataset is displayed and that the full simulation on this site covers 6576 hours.} \texttt{2S-SP} and \texttt{2S-SP-C} have the closest trajectories to \texttt{P-MPC}, displayed in light gray. The gains obtained using the challengers over \texttt{HEU} are achieved in two different ways:
\paragraph{Limiting the number of subscribed power limit overruns in high load conditions} Most of the savings are obtained in the very first (0-25) and last (>200) hours of the simulation slice shown in Fig. \ref{fig:trajectories}. The high net load forces all algorithms to exceed $\overline{E}$ to meet the energy demand, even when not charging the battery. This penalty being fixed, charging the battery on top of importing the energy necessary to satisfy the demand that already exceeds the subscribed power limit does not generate any additional cost (apart from the cost due to a larger volume being bought). The grid import trajectory under the policy of \texttt{HEU} constantly stays over $\overline{E}$ during those two time-frames, while it is more often below when following the policies of the challengers. On one hand, \texttt{HEU} constantly pays the penalty for exceeding $\overline{E}$ and imports the electricity necessary to meet the demand without charging the battery. On the other hand, the challengers repeatedly charge and discharge the battery, alternating between time-steps when the battery is charging and when the penalty is paid, and time-steps when discharging the battery allows the grid import to stay below $\overline{E}$. In particular, notice how \texttt{2S-SP-C} succeeds in keeping the grid import below $\overline{E}$ in the 205-210 time frame.
\paragraph{Arbitrage on the electricity prices and high net load anticipation}  On the rest of the slice, additional savings are obtained thanks to the arbitrage done on the electricity prices. Indeed, they typically charge the battery during off-peak hours, when the price is low, and when the demand is low enough for the grid import (increased by the energy charged into the battery) to stay below $\overline{E}$. The battery is then discharged when the electricity prices are higher, or when it allows the grid import that would otherwise exceed $\overline{E}$ power to stay below, by satisfying a bit of this demand using the energy stored in the battery. In the first case, this saves a cost equal to the amount of energy discharged, multiplied by the price spread. In the second case, this spares the penalty cost associated with the overrun of the subscribed power limit.

\section{Conclusion}
We introduced four different multistage stochastic control algorithms, namely, model predictive control, stochastic programming, two-stage stochastic programming, and two-stage stochastic programming with scenario clustering. We proposed a generic mathematical formulation of those algorithms under the shared framework of Multistage Stochastic Programming. We presented a simple microgrid model and benchmarked our four challengers against two baselines: a simple heuristic, and a model predictive control using perfect forecasts. The benchmark is done through 61 simulations, each having a different microgrid setup and exogenous uncertainty distribution. Our two-stage stochastic programming algorithm shows similar performances to our stochastic programming one, while at the same time being much easier to implement (as it doesn't require any scenario tree reduction step) and less computationally heavy. Both algorithms outperform the model predictive control approach, reducing the electricity cost by 2.2 \% (which amounts to 7413\euro~of average yearly cost savings) when used over model predictive control, and only generate an extra cost of 1.5 \% compared to the perfect baseline.
We believe this is an important finding, which questions the financial interest in relying on third-parties services, such as weather/irradiance forecasts, to improve the net load forecast accuracy.

\bibliographystyle{unsrt}  
\bibliography{references}

\end{document}